\documentclass{article}
\usepackage{amsthm}
\usepackage{amsfonts}

\begin{document}

Tsemo Aristide

PKFOKAM Institute OF EXCELLENCE

P.O. BOX 11646, Yaounde, Cameroon.

tsemo58@yahoo.ca

\bigskip
\bigskip

\centerline{\bf Automorphisms of projective manifolds.}

\bigskip

\centerline{\bf Abstract.}

{\it Let $(M,P\nabla_M)$ be a compact projective manifold and $Aut(M,P\nabla_M)$ its group of automorphisms. The purpose of this paper is to study the topological properties of $(M,P\nabla_M)$ if $Aut(M,P\nabla_M))$ is not discrete by applying the results that I have shown in [13] and the Benzekri's functor which associates to a projective manifold a radiant affine manifold. This enables us to show that the orbits of the connected component of $Aut(M,P\nabla_M)$ are immersed projective submanifolds. We also classify $3$-dimensional compact projective manifolds such that $dim(Aut(M,P\nabla_M))\geq 2$.  }

\bigskip
\bigskip

\centerline{\bf 1. Introduction}

\bigskip

  The purpose of this paper is to study the group of automorphisms of projective manifolds. Firstly we recall the definition of $(X,G)$ manifolds, their group of automorphisms and morphisms between $(X,G)$-structures. We applied the results described in the general framework of $(X,G)$-manifolds to the category of affine manifolds and projective manifolds. Benzekri has constructed a functor which associates to a projective manifold $(M,P\nabla_M)$ a radiant affine manifold $(B(M),\nabla_{B(M)})$ whose underlying topological space is $M\times S^1$. It enables us to show that there exists a surjective morphism between the connected component $Aut(B(M),\nabla_{B(M)})_0$ of the group of affine automorphisms of $(B(M),\nabla_{B(M)})$ and the connected component $Aut(M,P\nabla_M)_0$ of the group of projective automorphisms of $(M,P\nabla_M)$.
  
   Let $(M,\nabla_M)$ be a compact affine manifold, in [13], I have studied the relations between $Aut(M,\nabla_M)$  and the topology of $M$. This enables us to show that the orbits of $Aut(M,P\nabla_M)_0$ are projective immersed submanifolds. In the last section, we study the automorphisms group of $2$ and $3$ dimensional projective manifolds. We remark that a $2$-dimensional projective manifold whose group of automorphisms is not discrete is homeomorphic to the sphere, the $2$-dimensional projective space or the two dimensional torus. Finally we show that a $3$-dimensional projective manifold $(M,P\nabla_M)$ whose developing map is injective and such that $dim(Aut(M,P\nabla_M)\geq 2$ is homeomorphic to a spherical manifold,  $S^2\times S^1$, or a finite cover of $M$ is the total space of a torus bundle.
   
 Remark that $(X,G)$ manifolds play an important role in low dimensional topology: seven of the eight geometry of Thurston are examples of projective geometry (see Cooper and Goldman [8] p. 1220). In [12] p.17, Sullivan and Thurston note that the existence of a $(X,G)$-structure on every $3$-manifold implies the Poincare conjecture.

\bigskip

{\bf 2. $(X,G)$-manifolds.}

\bigskip

A  $(X,G)$ model is a finite dimensional differentiable manifold  $X$, endowed with an effective and transitive action of a Lie group $G$   which satisfies the unique extension property. This is equivalent to saying that:  two elements $g,g'$ of  $G$ are equal if and only if  their respective  restriction   to a non empty open subset  of $X$ are equal. 

A $(X,G)$ manifold $(M,X,G)$ is a differentiable manifold $M$, endowed with an open covering $(U_i)_{i\in I}$ such that for every $i\in I$, there exists a differentiable map $f_i:U_i\rightarrow X$ which is a diffeomorphism onto its image and  $f_i\circ f_j^{-1}$ coincides with the restriction of an element $g_{ij}$ of $G$ to $f_j(U_i\cap U_j)$. The map $f_i$ is called a $(X,G)$ chart.

A $(X,G)$ structure defined on $M$ can be lifted to the universal cover  $\tilde M$ of $ M$. This structure  is defined by a local diffeomorphism $D_M:\tilde M\rightarrow X$. This implies that a $(X,G)$ chart of  this structure is  an open subset $U$ of $\tilde M$ such that the restriction of $D_M$ to $U$ is a diffeomorphism onto its image.

\medskip

Let $(X,G)$ and $(X',G')$ be two models, $\phi:X\rightarrow X'$ a differentiable map and  $\Phi:G\rightarrow G'$ a morphism of groups such that for every $g\in G$, the following diagram is commutative:

$$
\matrix{ X&{\buildrel{g}\over{\longrightarrow}}& X\cr \phi\downarrow &&\downarrow \phi\cr X'&{\buildrel{\Phi(g)}\over{\longrightarrow}}& X'}
$$

Let $(M,X,G)$ (resp. $(M',X',G')$) be  an $(X,G)$ manifold (resp. an $(X',G')$ manifold). A $(\Phi,\phi)$-morphism $f:(M,X,G)\rightarrow (M',X',G')$ is a differentiable map: $f:M\rightarrow M'$ such that for every chart $(U_i,f_i)$ of $M$ such that $f(U_i)$ is contained in the chart $(V_j,f'_j)$ of $M'$, there exists an element $g\in G$ such that the restrictions of $f'_j\circ f\circ f_i^{-1}$ and $\Phi(g)\circ\phi $  to $f_i(U_i)$ coincide.

We will denote by $Aut(X,M,G)$ the group of $(Id_G,Id_X)$-automorphisms of $(M,X,G)$ and by $Aut(M,X,G)_0$ its connected component. It is a Lie group endowed with the compact open topology. For every element $g\in Aut(\tilde M,X,G)$, the developing map defines a representation $H_M:Aut(\tilde M,X,G)\rightarrow G$ such that the following diagram is commutative:

$$
\matrix{ \tilde{M}&{\buildrel{g}\over{\longrightarrow}}& \tilde M\cr D_M\downarrow &&\downarrow D_M\cr X&{\buildrel{H_M(g)}\over{\longrightarrow}}& X}
$$

Remark that the group of Deck transformations that we identify to the fundamental group $\pi_1(M)$, of $M$, is a subgroup of $Aut(\tilde M,X,G)$. The restriction $h_M$ of $H_M$ to the fundamental group $\pi_1(M)$, of $M$ is called the holonomy representation of the $(X,G)$ manifold $(M,X,G)$.

\bigskip

The pullback $p_M(f)$ of an element $f$ of $Aut(M,X,G)$, by the universal covering map, $p_M:\tilde M\rightarrow M$ is an element of $Aut(\tilde M,X,G)$  which belongs to the normalizer $N(\pi_1(M))$ of $\pi_1(M)$ in $Aut(\tilde M,X,G)$. Conversely, every element $g$ of $N(\pi_1(M))$ induces an element $A_M(g)$ of $Aut(M,X,G)$ such that the following diagram is commutative:

$$
\matrix{ \tilde{M}&{\buildrel{g}\over{\longrightarrow}}& \tilde M\cr p_M\downarrow &&\downarrow p_M\cr M&{\buildrel{A_M(g)}\over{\longrightarrow}}& M}.
$$

The kernel of the morphism $A_M:N(\pi_1(M))\rightarrow Aut(M,X,G)$ is $\pi_1(M)$ and $A_M$ is a local diffeomorphim. We will denote  by  $N(\pi_1(M))_0$ the connected component of $N(\pi_1(M))$, it is also the connected component of the commutator of $\pi_1(M)$ in $Aut(\tilde M,X,G)$. Since $A_M$ is locally invertible, it induces an isomorphism between the Lie algebra $n(\pi_1(M))$ of $N(\pi_1(M))$ and the Lie algebra $aut(M,X,G)$ of $Aut(M,X,G)$.
If  $(M,X,G)$ is a compact $(X,G)$ manifold, $aut(M,X,G)$ is isomorphic to space of elements of ${\cal G}$, the Lie algebra of $G$, which are invariant by $h_M(\pi_1(M))$.

\bigskip

{\bf 3. Affine and projective structures.}

\bigskip

Let $\mathbb{R}^n$ be the $n$-dimensional real vector space. We denote by $Gl(n,\mathbb{R})$ the group of linear automorphisms of $\mathbb{R}^n$ and by $Aff(n,\mathbb{R})$ its group of affine transformations. If we fix an origin $0$ of $\mathbb{R}^n$,  for every element $f\in Aff(n,\mathbb{R})$, 
we can write $f=(L(f),a_f)$ where $L(f)$ is an element of $Gl(\mathbb{R}^n)$ and $a_f=f(0)$. The couple $(\mathbb{R}^n,Aff(n,\mathbb{R}))$ is a model. A $(\mathbb{R}^n, Aff(n,\mathbb{R}))$ manifold is also called an affine manifold.
Equivalently, an $(\mathbb{R}^n,Aff(n,\mathbb{R}))$ manifold is a $n$-dimensional  differentiable manifold $M$ endowed with a connection $\nabla_M$ whose curvature and torsion tensors vanish identically.

Remark that the linear part $L(h_M)$ of the holonomy representation $h_M$ of an affine manifold $(M,\nabla_M)$ is the holonomy of the connection $\nabla_M$. We say that the $n$-dimensional affine manifold $(M,\nabla_M)$ is radiant if its holonomy $h_M$ fixes an element of $\mathbb{R}^n$, this is equivalent to saying that $h_M$ and $L(h_M)$ are conjugated by a translation.

The $n$-dimensional real projective space $\mathbb{R}P^n$ is the quotient of $\mathbb{R}^{n+1}-\{0\}$ by the equivalence relation defined by $x\simeq y$ if and only there exists $\lambda\in\mathbb{R}$ such that $x=\lambda y$. If $x$ is an element of $\mathbb{R}^{n+1}-\{0\}$, 
we will denote by $[x]_{\mathbb{R}P^n}$ its equivalent class.  The group $Gl(n+1,\mathbb{R})$ acts transitively on $\mathbb{R}P^n$ by the action defined by $g.[x]_{\mathbb{R}P^n}=[g.x]_{\mathbb{R}P^n}$ the kernel of this action is the group $H_{n+1}$ of 
homothetic maps. We denote by $PGl(n+1,\mathbb{R})$ the quotient $Gl(n+1,\mathbb{R})$ by $H_n$. The couple $(\mathbb{R}P^n,PGl(n+1,\mathbb{R}))$ is a model.  A $(\mathbb{R}P^n,PGl(n+1,\mathbb{R}))$ is also called a projective manifold. Equivalently, a projective manifold can be defined by a differentiable manifold $M$ endowed with a projectively flat connection $P\nabla_M$. We will denote it by $(M,P\nabla_M)$.

The $n$-dimensional sphere $S^n$ is the quotient of $\mathbb{R}^{n+1}-\{0\}$ by the equivalence relation defined by $x\simeq y$ if and only if there exists $\lambda>0$ such that $x=\lambda y$. Let $x$ be an element of $\mathbb{R}^{n+1}-\{0\}$, we will denote by $[x]_{S^n}$ its equivalence class for this relation. Remark that if $\langle,\rangle$ is an Euclidean metric defined on $\mathbb{R}^{n+1}$, there exists a bijection between  the unit sphere $S^n_{\langle,\rangle}=\{x:x\in\mathbb{R}^{n+1}, \langle x,x\rangle=1\}$ and $S^n$ defined by the restriction of the equivalence relation to $S^n_{\langle,\rangle}$.

There exists a map $D_{S^n}:S^n\rightarrow \mathbb{R}P^n$ such that for every element $x$ of $\mathbb{R}^{n+1}-\{0\}$,  $[x]_{\mathbb{R}P^n}=D_{S^n}([x]_{S^n}).$ The map $p_n$ is a covering, thus is the developing map of    a projectively flat connection $P\nabla_{S^n}$ defined on $S^n$.

\medskip

A $p$-dimensional  projective submanifold $(F,P\nabla_F)$ of the projective manifold $(M,P\nabla_M)$ is a $p$-dimensional submanifold $F$ of $M$ endowed with a structure of a projective manifold, such that the canonical embedding $i_F:(F,P\nabla_F)\rightarrow (M,P\nabla_M)$ is a morphism of projective manifolds.

\medskip

Let $\hat F$ be the universal cover of $F$, we can lift $i_F$ to a projective map $\hat i_F:\hat F\rightarrow \hat M$. The image of $D_M\circ \hat i_F$ is contained in a $p$-dimensional projective subspace $U_F$ of $\mathbb{R}P^n$. The map $D_M\circ \hat i_F:\hat F\rightarrow U_F$ is a developing map of $F$.
There exists a canonical  morphism $\pi_F:\pi_1(F)\rightarrow \pi_1(M)$ induced by $i_F$. Let $\gamma$ be an element of $\pi_1(F)$, the holonomy $h_F(\gamma)$ is the restriction of $h_M(\pi_F(\gamma))$ to $U_F$. If there is no confusion, we are going to denote $h_M(\pi_F(\gamma))$ by $h_M(\gamma)$.

\medskip

{\bf Proposition 3.1.}
{\it The group of automorphisms of the $n$-dimensional projective manifold $S^n$ is isomorphic to  $Sl(n+1,\mathbb{R})$, the group of invertible $(n+1)\times (n+1)$ matrices such that for every element $A\in Sl(n+1,\mathbb{R}), |det(A)|=1$.}

\medskip

{\bf Proof.}
 Let $g$ be an element of $Sl(n+1,\mathbb{R})$. For every $[x]_{S^n}\in S^n$, we write $u_g(x)=[g(x)]_{S^n}$. Let $[g]$ be the image of $g$ by the quotient map $Sl(n+1,\mathbb{R})\rightarrow  PGl(n+1,\mathbb{R})$, we have $[g]\circ D_{S^n}=D_{S^n}\circ u_g$. This implies that $u_g$ is an element of $Aut(S_n,P\nabla_{S^n})$. Suppose that $u_g=Id_{S^n}$, it 
 implies that for every $[x]\in S^n$,  $g(x)=\lambda(x)x,\lambda(x) >0$,  we deduce that $g(x)=\lambda Id_{\mathbb{R}^n},\lambda>0$,
  and $\lambda^{n+1}=1$ since $g\in Sl(n+1,\mathbb{R})$. This implies that $\lambda =1$.
  We deduce that $u:Sl(n+1,\mathbb{R})\rightarrow Aut(S^n,P\nabla_{S^n})$ defined by $u(g)=u_g$ is injective.
Let $f$ be an element of $Aut(S^n,P\nabla_{S^n})$, there exists an element $[g]\in PGl(n+1,\mathbb{R})$ such that $[g]\circ D_{S^n}=D_{S^n}\circ f$. Consider an element $g\in Sl(n+1,\mathbb{R})$ whose image by the quotient map is $[g]$, $f=u_g$. This implies that $u$ is an isomorphism.

\bigskip

{\bf The Benzecri correspondence.}

\medskip

Consider the embedding $i_n^G:Gl(n,\mathbb{R})\rightarrow Gl(n+1,\mathbb{R})$ defined by $i_n^G(A)=\pmatrix{A & 0\cr 0 &1}$ and the open embedding $i_n:\mathbb{R}^n\rightarrow \mathbb{R}P^n$ defined by $i_n(x_1,...,x_n)=[x_1,...,x_n,1]$. For every elements $x\in\mathbb{R}^n$ and $g\in Gl(n,\mathbb{R})$, we have $i_n(g(x))=i_n^G(g)(i_n(x))$. We deduce that for every affine manifold  $(M,\nabla_M)$  whose developing map is $D_M$, there exists a projective structure defined on  $M$ whose developing map is $i_n\circ D_M$. 

Benzecri [4] p.241-242 has defined a functor between the category of projective manifolds of dimension $n$ and the category of radiant affine manifolds of dimension $n+1$ which can be described as follows:

Firstly, we remark that since the universal cover $\tilde M$ of the projective manifold $M$ is simply connected and $p_n:S^n\rightarrow P\mathbb{R}^n$ is a covering map, the theorem 4.1 of Bredon [5]  p.143 implies that the development map $D_M:\tilde M\rightarrow P\mathbb{R}^n$, can be lifted to a local diffeomorphism  $D'_M:\tilde M\rightarrow S^n$ which is a projective morphism. Let $N(\pi_1(M))$ be the normalizer of $\pi_1(M)$ in $Aut(\tilde M,P\nabla_{\tilde M})$, for every $g\in N(\pi_1(M))$, there exists $H'_M(g)\in Aut(S^n,P\nabla_{S^n})$ such that the following diagram is commutative: 

$$
\matrix{ \tilde{M}&{\buildrel{g}\over{\longrightarrow}}& \tilde M\cr D'_M\downarrow &&\downarrow D'_M\cr S^n&{\buildrel{H'_M(g)}\over{\longrightarrow}}& S^n}.
$$

We will denote by $h'_M$ the restriction of $H'_M$ to $\pi_1(M)$.

\medskip

There exists a local diffeomorphism $D_{\tilde M\times S^1}:\tilde M\times \mathbb{R}^*_+\rightarrow \mathbb{R}^{n+1}-\{0\}$ defined by $D_{\tilde M\times S^1}(x,t)=tD'_M(x)$, which is  the developing map of a radiant structure defined on $M\times S^1$ whose holonomy representation $h_{M\times S^1}:\pi_1(M\times S^1)\rightarrow Gl(n+1,\mathbb{R})$ is defined by $h_{M\times S^1}(\gamma,n)=2^nh'_M(\gamma)$. This radiant affine manifold $M\times S^1$ is the construction of Benzecri, we will often denote this affine structure by $(B(M),\nabla_{B(M)})$ and by $p_{B(M)}:M\times S^1\rightarrow M$ the projection on the first factor.

Let $f:(M,P\nabla_M)\rightarrow (N,P\nabla_N)$ be a morphism between $n$-dimensional projective manifolds; $f$ can lifted to the projective  the morphism $\tilde f:\tilde M\rightarrow \tilde N$. We deduce the existence of a morphism of affine manifolds $f':\tilde M\times \mathbb{R}_+^*\rightarrow \tilde N\times \mathbb{R}_+^*$ defined by $f'(x,t)=(\tilde f, t)$. The morphism $f'$ is equivariant with respect to the action of $\pi_1(B(M))$ on $\tilde M\times\mathbb{R}^*_+$ and $\pi_1(B(N))$ on $\tilde N\times\mathbb{R}^*_+$, and covers a morphism $b(f):B(M)\rightarrow B(N)$.

Let $(N,\nabla_N)$ be a $n$-dimensional  radiant affine manifold. We suppose that the holonomy of $N$ fixes the origin of $\mathbb{R}^n$. The vector field defined on $\mathbb{R}^n$ by $X_R^{\mathbb{R}^n}(x)=x$  is invariant by the holonomy. Its pullback by the developing map is a vector field $X_R^{\tilde N}$ of $\tilde N$ invariant by $\pi_1(N)$. We deduce that $X_R^{\tilde N}$ is the pullback of a vector field $X_R^N$ of $N$ called the radiant vector field of $N$.

 \bigskip
 
{\bf Proposition 3.2.}
{\it Let $(M,P\nabla_M)$ be a compact projective manifold. There exists a surjective morphism of groups between the connected component of $Aut(B(M),\nabla_{B(M)})$ and the connected component of $Aut(M,P\nabla_M)$. }

\medskip

{\bf Proof.}
Let $f$ be an element of $Aut(B(M),\nabla_{B(M)})_0$, the connected component of $Aut(B(M),\nabla_{B(M)})$. 
 Consider an element $\tilde f$ of $Aut(\tilde M\times\mathbb{R}^*_+)_0$ over $f$. For every $\tilde x\in\tilde M$ and $t\in \mathbb{R}^*_+$, we can write $\tilde f(\tilde x,t)=(\tilde g(\tilde x,t), h(\tilde x,t))$. The flow of $X_R^{\mathbb{R}^{n+1}}$ is in the center of $Gl(n+1,\mathbb{R})$, we deduce that $\tilde f $ commutes with the flow $X_R^{\tilde B(M)}$, $g(\tilde x,t)$ does not depend of $t$ and $h(\tilde x,t)=th(\tilde x,1)$.

Let $\gamma$ be an element of $\pi_1(M)$, since $(\gamma,2).(\tilde x,t)=(\gamma(\tilde x), 2t)$ is an element of $\pi_1(B(M))$ and $\tilde f$ is an element of $N(\pi_1(B(M)))_0$, we deduce that $(\gamma,2)$ commutes with $\tilde f$ and $\tilde g$ commute with $\gamma$. This implies that there exists an element $g$ of $Aut(M,P\nabla_M)$ whose lifts is $\tilde g$. Remark that since $\tilde f$ is an affine transformation,  $h(x,1)$ is a constant. The correspondence $P:Aut(B(M),\nabla_{B(M)})_0\rightarrow Aut(M,P\nabla_M)_0$ defined by $P(f)=g$ is well defined and is surjective morphism of groups since for every element $f\in Aut(M,P\nabla_M)_0$, $P(b(f))=f$.

\bigskip
 Let $(M,P\nabla_M)$ be a  projective manifold $M$, the orbits of the radiant flow $\phi_t^{B(M)}$ of $X_R^{B(M)}$  are compact. The images of the elements of $\phi_t^{B(M)}$ by $P$  are the identity on  $(M,P\nabla_M)$. This implies that $dim(Aut(M,P\nabla_M))+1\leq dim(Aut(B(M),\nabla_{B(M)})))$.
 We deduce that if $(M,P\nabla_M)$ is a projective manifold, such  that $Aut(M,P\nabla_M)$ is not discrete,  the dimension of $Aut(B(M),\nabla_{B(M)})$ is superior or equal to $2$.

\bigskip
{\bf 4. Automorphisms of projective manifolds and automorphisms of radiant affine manifolds.}

\bigskip

Let $(N,\nabla_N)$ be an affine manifold. In [13], I have shown that  $aut(N,\nabla_N)$, the Lie algebra of $Aut(N,\nabla_N)$ is endowed with an associative product defined by $X.Y={\nabla_M}_XY$. We deduce that  $\nabla_{B(M)}$ defines on $aut(B(M),\nabla_{B(M)})$ 
an associative structure which can be pulled back to $n(\pi_1(B(M))$. It results that the Lie algebra $H_{B(M)}(n(B(M),\nabla_{B(M)}))$ of the image of $N(\pi_1(B(M)))$ by $H_{B(M)}$ is stable by the canonical product of matrices which is the image of the associative product of $n(\pi_1(B(M)))$ by $H_{B(M)}$. Remark that $H_{B(M)}(n(B(M),\nabla_{B(M)}))$ is isomorphic to $n(B(M),\nabla_{B(M)})$. The theorem 23 of chap. III of [1] implies that we can write: $H_{B(M)}(n(B(M),\nabla_{B(M)}))=S_M\oplus N_M$ where $S_M$ is a semi-simple associative algebra and $N_M$ a nilpotent associative algebra.

In [14],  by using this associative product, I have shown that the orbits of the canonical action of  $Aff(N,\nabla_N)_0$ on $N$ are immersed affine submanifolds of $(N,\nabla_N)$ and are  the leaves of a (singular) foliation. This leads to the following result:

\medskip

{\bf Proposition 4.1.}
{\it Let $(M,P\nabla_M)$ be a projective manifold. The orbits of the  action of $Aut(M,P\nabla_M)_0$ on $M$ are immersed projective submanifolds and are the  leaves of a singular foliation.}

\medskip

{\bf Proof.}
The orbits of $Aut(B(M),\nabla_{B(M)})_0$ are immersed affine submanifolds of $B(M)$. The proposition 3.2 shows that there exists a surjective map $P:Aut(B(M),\nabla_{B(M)})_0\rightarrow Aut(M,P\nabla_M)_0$ such that, for every $g\in Aut(B(M),\nabla_{B(M)})_0$ and $x\in B(M)$, $p_{B(M)}(g(x))=P(g)(p_B(x))$. This implies that the orbits of $Aut(M,P\nabla_M)_0$ are the images of the orbits of $Aut(B(M),\nabla_{B(M)})_0$ by the quotient map $B(M)\rightarrow M$.

\bigskip 

{\bf Theorem 4.1}
{\it Let $(M,P\nabla_M)$ be a compact oriented projective manifold  of dimension superior or equal to $2$. Suppose that  $H_M(N(\pi_1(M)))$ acts transitively on $\mathbb{R}P^n$, then $(M,P\nabla_M)$ is isomorphic to a finite quotient of $\mathbb{K}P^m$ by a subgroup of $\mathbb{K}$ where $\mathbb{K}$ is the field of real numbers, complex numbers, quaternions or octonions. The action of $\pi_1(M)$ on $\mathbb{K}P^n$ is induced by its action on $\mathbb{K}^{m+1}$ by homothetic maps.}

\medskip

{\bf Proof.}
The fact that $H_M(N(\pi_1(M)))$ acts transitively on $\mathbb{R}P^n$ implies that $H'_M(N(\pi_1(M)))$ acts transitively on $S^n$. The theorem of Montgomery Zipplin [11] p.226 implies that a connected compact subgroup $K'$ of $H'_M(N(\pi_1(M)))$ acts transitively on $S^n$. The theorem I p. 456 of Montgomery and Samelson [10] implies that a connected compact simple subgroup $C'$ of $K'$ acts transitively on $S^n$. The Lie algebra of the connected component $C$ of ${H'}_M^{-1}(C')$ is isomorphic to the Lie algebra of $C'$ since the kernel of $H'_M$ is discrete. This implies that $C$ is compact. Remark that the orbits of the action of $C$ on $\tilde M$ are open. We deduce that $C$ acts transitively on $\tilde M$ and $\tilde M$ is compact. This implies that $D'_M:\tilde M\rightarrow S^n$ is a covering since it is a local diffeomorphism defined between compact manifolds. This  implies that $D'_M$ is a diffeomorphim since $S^n$ and $\tilde M$ are simply connected.   

We can write $\mathbb{R}^n=\oplus_{i\in I}U_i$ where $U_i$ is an irreducible component of the action of $h'_M(\pi_1(M))$ since $\pi_1(M)$ is finite. 

Let $i,j\in I$, consider two non zero elements $x_i\in U_i, x_j\in U_j$, since $H'_M(N(\pi_1(M))_0)$ acts transitively on $S^n$ there exists $B\in H'_M(N(\pi_1(M))_0)$ such that $B(x_i)=cx_j$, $B(U_i)\cap U_j$ is invariant by $H'_M(\pi_1(M))$, we deduce that $B(U_i)=U_j$ since $U_j$ is irreducible and $B$ is an isomorphism.

The group of automorphisms of the irreducible representation $U_i$ is $\mathbb{K}$ where $\mathbb{K}=\mathbb{R}, \mathbb{C}, \mathbb{H}$ or $\mathbb{O}$. We deduce that $\mathbb{R}^{n+1}$ is a $\mathbb{K}$ vector space and the action of $H'_M(\pi_1(M))_0)$ on $\mathbb{K}^n$ is induced by its action on $\mathbb{K}$ by right multiplication of elements of $\mathbb{K}$.

\medskip

{\bf Remark.}

\medskip

Suppose that the dimension of $M$ is even,  and $H_M(N(\pi_1(M))_0))$ acts transitively on $\mathbb{R}P^n$. The proof of the previous theorem  can be simplified  as follows: Every element of $H'_M(\pi_1(M))$ has a fixed point since every element of $Gl(2n+1,\mathbb{R})$ has a real eigenvalue. We deduce that $H'_M(\pi_1(M))$ is the identity and there exists a map $f:M\rightarrow \mathbb{R}P^n$ such that $D_M=f\circ p_M$. This implies that $f$ is a covering map and $M$ is homeomorphic to $S^n$ or $\mathbb{R}P^n$.

\medskip

Let $(M,P\nabla_M)$ be a projective manifold, suppose that $Aut(M,P\nabla_M)_0$ is not solvable. This implies that $Aut(B(M),\nabla_{B(M)})_0$ and the connected component of the normalizer $N(\pi_1(B(M)))_0$ of $\pi_1(B(M))$ in $Aut(\widetilde{B(M)},\nabla_{\widetilde{B(M)}})$ are not solvable. We deduce that the image  of $N(\pi_1(B(M))$ by $H_{B(M)}$ contains  a subgroup $H_{S^1}$ isomorphic to $S^1$. We denote by $X'_{B(M)}$, a vector field which generates the Lie algebra of $H_{S^1}$, its pullback by  the developing map $D_{B(M)}$ of $B(M)$ is a vector field $\tilde X_{B(M)}$ invariant by the fundamental group of $B(M)$. We deduce that there exists a vector field $X_{B(M)}$ of $B(M)$ whose pullback by the universal covering map is $\tilde X_{B(M)}$. Suppose that the developing map is injective,  the flow of $X_{B(M)}$ defines an action of $S^1$ on $B(M)$ which is transverse  and commutes with  the radial flow. This implies there exists a vector $X_M$ on $M$ which is the image of $X_{B(M)}$ by the map induced by $p_{B(M)}:B(M)\rightarrow M$. The vector field $X_M$ induces an action of $S^1$ on $M$: We have:

\medskip

{\bf Proposition 4.2.}
{\it Let $(M,P\nabla_M)$ be a compact projective manifold whose developing map is injective, suppose that $Aut(M,P\nabla_M)_0$ is not solvable, then  $M$ is endowed with a non trivial action of $S^1$.}

\bigskip

{\bf 5. Automorphisms of projective manifolds of dimension $2$ and $3$. }

\bigskip

 In dimension $2$, we have the following result:

\medskip

{\bf Proposition 5.1.}
{\it Let $(M,P\nabla_M)$ be a $2$-dimensional compact connected oriented projective manifold, suppose that $Aut(M,P\nabla_P)$ is not discrete, then $M$ is homeomorphic to the $2$-dimensional torus or to the sphere.}

\medskip

{\bf Proof.}
Suppose $N_M\neq 0$, there exists a non zero element $A_M\in N_M$ such  that $A_M^2=0$, we deduce that $dim(ker(A_M))=2, dim(Im(A_M))=1$. Remark that $Im(A_M)$ is fixed by the holonomy.

Suppose that $N_M=0$, we deduce that $dim(S_M)\geq 2 $, there exists a non zero element distinct of the identity $e_M\in S_M$ such that $e^2_M=e_M$. To see this remark that $S_M$ contains either an associative algebra isomorphic to the associative algebra of $2\times 2$ real matrices or two idempotents which are linearly independent. The linear map $e_M$ is diagonalizable and its eigenvalues are equal to $0$ and $1$. Since the flow of $e_M$ is distinct of the radial flow, we deduce that $0$ is an eigenvalue of $e_M$. This implies that either the dimension of the eigenspace associated to $0$ is $1$, or the the dimension of the eigenspace associated to $1$ is $1$. We deduce that the holonomy preserves a vector subspace of dimension $1$. 

We conclude that if $Aut(M,P\nabla_M))$ is not discrete, its holonomy fixed a point of $P\mathbb{R}^2$. The lemma 2.5 p. 808 in Goldman [9] implies that the Euler number of $M$ is positive.

\bigskip

{\bf Dimension $3$.}

\medskip

In this section, we study the group of automorphisms of a connected $3$-dimensional compact projective manifold $(M,P\nabla_M)$ whose group of automorphisms is not discrete.

\medskip

{\bf $Aut(M,P\nabla_M)_0$ is not solvable.}

\medskip

Suppose that $Aut(M,P\nabla_M)_0$ is not solvable, then $Aut(B(M),\nabla_{B(M)})_0$ and $N(\pi_1(B(M)))_0$ are not solvable. We deduce that the connected  subgroup of $Gl(n+1,\mathbb{R})$, $H_{B(M)}(N(\pi_1(M))_0)$ contains a subgroup $H"$ isomorphic to $S^1$.  We denote by $X"_{B(M)}$ a vector field which generates the Lie algebra of $H"$. The pullback $X'_{B(M)}$ of $X"_{B(M)}$ by $D_{B(M)}$ is the pullback of a vector field $X_{B(M)}$ of $B(M)$ by $p_{B(M)}$.

\medskip

Suppose that the set of fixed points of $H"$ is not empty, we can write $\mathbb{R}^4=U\oplus V$ where $U$ is a $2$-dimensional vector subspace corresponding to the non trivial irreducible submodule of $H"$ and $V$ the set of fixed points. Remark that $h_{B(M)}(\pi_1(B(M)))$ preserves $U$ and $V$ since it commutes with $H"$. This implies that there exists a foliation ${\cal F}_U$ (resp. ${\cal F}_V$) on $B(M)$ whose pullback by the universal covering map is the pullback  by $D_{B(M)}$ of the foliation of $\mathbb{R}^4$ whose leaves are $2$-dimensional affine spaces parallel to $U$ (resp. parallel to $V$).

\medskip

{\bf Proposition 5.2.}
{\it Suppose that $V\cap D_{B(M)}(\widetilde{B(M)})$ is empty. Then a finite cover of $M$ is a total space of a fibre bundle over $S^1$ whose fibre is $T^2$.}

\medskip

{\bf Proof.}
The vector field defined by $Y"(u,v)=u; u\in U, v\in V$ is invariant by the holonomy of $B(M)$. To show this, remark that the restriction of $H"$ to $U$ defines on it a complex structure and since $h_{B(M)}(\pi_1(B(M)))$ commutes with $H"$, its restriction to $U$ are morphisms of that complex structure. The pullback of $Y"$    by $D_{B(M)}$ is the pullback of a vector field $Y_{B(M)}$ of $B(M)$ by the universal covering map. The image $Y_M$ of $Y_{B(M)}$ by $p_{B(M)}$ and $X_M$ commute and generate a locally free action of $\mathbb{R}^2$ on $M$ since $V\cap D_{B(M)}(\widetilde{B(M)})$ is empty.
 Chatelet, Rosenberg and Weil [6] implies that $M$ is the total space of a fibre bundle over $S^1$ whose fibre is $T^2$.

\medskip

If ${\cal F}_V$ has compact leaf, we have the following result:

\medskip

{\bf Proposition 5.3.}
{\it Let $(M,P\nabla_M)$ be a $3$-dimensional compact projective manifold whose developing map is injective. Suppose that $Aut(M,P\nabla_M)_0$ is not solvable and $V\cap D_{B(M)}(\widetilde{B(M)})$ is not empty. Then the holonomy of $(M,P\nabla_M)$ is solvable.}

\medskip

{\bf Proof.}
Let $\hat F_0$ be a connected component of $V\cap D_{B(M)}(\widetilde{B(M)})$ its image by the universal covering map is a compact leaf $F_0$ compact leaf of  ${\cal F}_V$ which is a $2$-dimension compact  affine manifold, we deduce that its fundamental group is solvable. Let $r$ be the restriction of $h(\pi_1(B(M))$ to $V$, since $h_{B(M)}(\pi_1(B(M))$ preserves $V$, we have an exact sequence:

$$
1\rightarrow Ker(r)\rightarrow h_{B(M)}(\pi_1(B(M))\rightarrow Im(r)\rightarrow 1
$$
 The groups $Ker(r)$ is solvable since it restriction to $U$ commutes with a non trivial linear action of $S^1$. The group $Im(r)$ is also solvable since it is contained in $h_{F_0}(\pi_1(F_0))$, we deduce that $h_{B(M)}(\pi_1(B(M))$ is solvable.

\bigskip

{\bf $Aut(M,P\nabla_M)_0$ is solvable.}

\medskip

In this section we study $3$-dimensional projective manifolds whose group of automorphisms is solvable. We can decompose the associative algebra $n(\pi_1(B(M))$ by writing:  $n(\pi_1(B(M))=S_M\oplus N_M$, where $S_M$ is a  semi-simple associative algebra and $N_M$ a nilpotent associative algebra. We deduce that $S_M$ is the direct product of associative algebras isomorphic to either $\mathbb{R}$ or $\mathbb{C}$ and is commutative. It results that the fact that $Aut(M,P\nabla_M)_0$ is not commutative implies that $N_M$ is not commutative.

\medskip

{\bf $Aut(M,P\nabla_M)_0$ is solvable and  is not commutative.}

\bigskip

{\bf Theorem 5.1.}
{\it Suppose that $N_M$ is not commutative, then  $h_{B(M)}(\pi_1(B(M)))$, the image of the holonomy of $B(M)$ is solvable.}

\medskip

{\bf Proof.}

First step:

Suppose that the square of every element of $N_M$ is zero.

Let $A,B\in N_M$ such that $AB\neq BA$.
 Suppose that $dim(ker(A))=3$. It implies that $dim(Im(A))=1$. Since $(A+B)^2=0$, we deduce that $AB+BA=0$ and $B(Ker(A))\subset Ker(A), B(Im(A))\subset Im(A))$, we deduce that the restriction of $B$ to $Im(A)$ is zero since $B$ is nilpotent and $dim(Im(A))=1$. This implies that $BA=0$, we deduce that $AB=0$ and $AB=BA$. Contradiction.

Suppose that that $dim(Ker(A))=dim(Im(A))=dim(Ker(B))=dim(Im(B))=2$. We deduce that $Im(A)=Ker(A), Im(B)=Ker(B)$ since $A^2=B^2=0$. If $Ker(A)\cap Ker(B)=0$, $\mathbb{R}^4= Ker(A)\oplus Ker(B)$ and $AB=BA=0$ since $AB+BA=0$. If $Ker(A)=Ker(B)$, $AB=BA=0$ since $Im(A)=Im(B)=Ker(A)=Ker(B)$. Contradiction.

We deduce that $dim(Ker(A)\cap Ker(B))=1$. We can write $\mathbb{R}^4=Vect(e_1,e_2,e_3,e_4)$ where $Vect(e_1)=Ker(A)\cap Ker(B)$, $Vect(e_1,e_2)=Ker(A)$ and $Vect(e_1,e_3)=Ker(B)$.  Every element in $h_{B(M)}(\pi_1(B(M))$ preserves $Ker(A)\cap Ker(B), Ker(A)$ and $Ker(A)+Ker(B)$ since it commutes with $A$ and $B$. We deduce that $\pi_1(B(M))$ is solvable since it preserves a flag.

Step 2.

Suppose that  there exists an element $A\in N_M$ such that $A^2\neq 0$.

 If $dim(Ker A)=1$, we have $(A^2)^2=0$ implies that $Im(A^2)\subset Ker(A^2)$. Remark that $x\in Ker(A^2)$ if and only if $A(x)\in Ker(A)$ and $x\in A^{-1}(Ker(A))$; $dim(A^{-1}(Ker(A))=2$ since $dim(Ker(A))=1$. We deduce that $dim(Ker(A^2))=dim(Im(A^2))=2$, and  $Ker(A)\subset Ker(A^2)=Im(A^2)\subset Im(A)$ and $\pi_1(B(M))$ preserves a flag since it commutes with $A$ and $dim(Im(A))=3$. 

Suppose that $dim (ker(A))=2$, we deduce that $dim(Im(A))=2$. Suppose that $Ker(A)\cap Im(A)=0$, we deduce that
$\mathbb{R}^4=Ker(A)\oplus Im(A)$. This is impossible since $A$ is nilpotent. We also remark that $Ker(A)$ is distinct of $Im(A)$ since $A^2\neq 0$. This implies that   $dim(Ker(A)\cap Im(A))=1$.
Every element of $\pi_1(B(M))$ preserves, $Ker(A)\cap Im(A), Ker(A)$ and $Ker(A)+Im(A)$ and thus preserves a flag. 
We deduce that $h_{B(M)}(\pi(B(M)))$ is solvable.

\bigskip

{\bf $Aut(M,P\nabla_M)_0$ is commutative.}

\bigskip

Suppose that the developing map is injective and $Aut(M,P\nabla_M)_0$ is commutative and its dimension is superior or equal to $2$. Let $X_M$ and $Y_M$ two  projective vector fields linearly independent. We denote by $X_{B(M)}$ and $Y_{B(M)}$ two  affine  vector fields of $B(M)$ whose respective images by $p_{B(M)}$ are $X_M$ and $Y_M$. There exist affine vector fields $X'_{B(M)}$ and $Y'_{B(M)}$ of $\mathbb{R}^4$ whose respective images by the covering map are $X_{B(M)}$ and $Y_{B(M)}$. Remark that if the group generated by $X_M$ and $Y_M$ acts freely on $M$, then Chatelet, Rosenberg and Weil [6] implies that $M$ is the total space of a torus bundle. 

In the rest of this section we assume that the set of zero of $X_M$ is not empty. This implies that we can assume that the set of zero $U$ of $X'_{B(M)}$ is not empty by eventually replacing $X'_{B(M)}$ with $X'_{B(M)}+cX_R$ where $c\in \mathbb{R}$ and $X_R$ is the radiant flow. We denote by $B(N)$ the image of $U\cap D_{B(M)}(\tilde M)$ by the covering map. Remark that $B(N)$ is not empty.

\medskip

{\bf Proposition 5.4.}
{\it Suppose that $dim(U)=3$ then $\pi_1(M)$ is solvable.}

\medskip

{\bf Proof.}
 Suppose that the restriction of $Y'_{B(M)}$ to $U$ is not zero. This implies that the restriction of $Y_{B(M)}$ to $B(N)$ is not zero. This implies that the group of projective automorphisms of $N$, the quotient of $B(N)$ by the radiant flow is not discrete. The proposition 5.1 implies that  $N$ has a finite cover  homeomorphic to $S^2$ or $T^2$. This implies that $\pi_1(B(N))$ is solvable. The restriction of $\pi_1({B(M)})$ to $U$ induces an exact sequence whose image is contained in the image of the holonomy representation of $B(N)$  and whose kernel is solvable. We deduce that $\pi_1(B(M))$ and $\pi_1(M)$ are solvable. 
 If the restriction of $Y'_{B(M)}$ to $U$ vanishes, let $V$ be the image of $X'_{B(M)}$, if $V$ is not contained in $U$, then $\mathbb{R}^4=U\oplus V$, and since $Y'_{B(M)}$ commutes with $X'_{B(M)}$, it preserves $V$. This implies that $X'_{B(M)}$ and $Y'_{B(M)}$ are linearly dependent contradiction.
 
 Suppose that $V$ is a subset of $U$, let $W$ be the image of $Y'_{B(M)}$, if $W$ is not contained in $U$, then we can apply the previous argument to obtain a contradiction by replacing $X'_{B(M)}$ by $Y'_{B(M)}$. Suppose that $W$ is contained in $U$, $V\cap W=\{0\}$ since $X'_{B(M)}$ and $Y'_{B(M)}$ are linearly independent. We deduce that, the holonomy of $B(M)$ preserves, $V, V\oplus W$ and $U$. This implies that the holonomy of $B(M)$ and $\pi_1(B(M))$ are solvable.

\medskip

{\bf Proposition 5.5.}
{\it Suppose that $dim(U)=2$, then $\pi_1(B(M))$ is solvable.}

\medskip

{\bf Proof.}
Suppose that $U\oplus Im(X'_{B(M)})=\mathbb{R}^4$.

Step $1$.
 
If the restriction of $X'_{B(M)}$ or $Y'_{B(M)}$to $Im(X'_{B(M)})$ are not a multiple of the identity, we deduce that $r(\pi_1(B(M))$, the image of  the restriction of $\pi_1(B(M))$ to $Im(X'_{B(M)})$ is solvable since it commutes with $X'_{B(M)}$ and $Y'_{B(M)}$. The restriction of $\pi_1(B(M))$ to $Ker(X_{B(M)})$ is contained in the holonomy group the $2$-dimensional closed affine manifold $B(N)$, we deduce that it is solvable. This implies $\pi_1(B(M))$ is solvable.

Step $2$.

Suppose that the restriction of $X'_{B(M)}$ and $Y'_{B(M)}$ to $Im(X'_{B(M)})$ are multiple of the identity. If the restriction of $Y'_{B(M)}$ to $U$ is equal to $aI_U$, we deduce that $Y'_{B(M)}$ is contained in the vector space generated by $Id_{\mathbb{R}^4}$ and $X'_{B(M)}$. This implies that the dimension of  the vector space generated by $X_M$ and $Y_M$  is $1$. Contradiction. We deduce that the restriction of $Y'_{B(M)}$ to $U$ is not a multiple of the $Id_U$. There exists a real $a$ such that the restriction of $Z=Y'_{B(M)}+aId_{\mathbb{R}^4}$ to $Im(X'_{B(M)})$ is zero. The restriction of $Z$ to $U$ is distinct of a multiple of the identity, we conclude that $\pi_1(B(M))$ is solvable by replacing $Y'_{B(M)}$ by $Z$ in the first step of the proof.

\medskip

Suppose that $dim(U\cap Im(X'_{B(M)})=1$. The vector subspaces, $U\cap Im(X'_{B(M)}), U, U\oplus Im(X'_{B(M)})$ are stable by the holonomy. We deduce that $\pi_1(B(M))$ preserves a flag and is solvable.

Suppose that $U=Im(X'_{B(M)})$.

We can write $\mathbb{R}^4=Vect(e_1,e_2,e_3,e_4)$ where $Vect(e_1,e_2)=U$ and $X'_{B(M)}(e_3)=e_1, X'_{B(M)}(e_4)=e_2$. Let $\gamma$ be an element of $\pi_1(B(M))$, if we write the fact that the matrix $M(\gamma)$ of $\gamma$ commutes with the matrix of $X'_{B(M)}$ in the basis $(e_1,e_2,e_3,e_4)$, we obtain that:

$$
M(\gamma)=\pmatrix{a_1 & b_1 & c_1 & d_1\cr a_2 & b_2 & c_2 & d_2\cr 0 & 0 & a_1 & b_1\cr 0 & 0 & a_2 & b_2}
$$

Since the restriction of $\pi_1(B(M))$ to $U$ is contained in the holonomy of $\pi_1(B(N))$ which is solvable, we deduce that  $\pi_1(B(M))$ is solvable.

\medskip

{\bf Proposition 5.6.}
{\it Suppose that $dim(U)=1$, then $\pi_1(B(M))$ is nilpotent.}

\medskip
{\bf Proof.}
We can write $n(\pi_1(B(M)) =S_M\oplus N_M$ where $S_M$ is semi-simple and $N_M$ nilpotent. Suppose that $N_M$ is not zero, and consider $n\in N_M$, if $n^2\neq 0$, $(n^2)^2=0$, $dim(Ker(n^2))\geq 2$, if $n^2=0$, $dim(Ker(n))\geq 2$. We can apply the proposition 5.4 and the proposition 5.5.   

Suppose that $N_M=0$,  $n(B(M))$ contains a non zero idempotent $u_1$ which is not a multiple of the identity, since it does not generates the radiant flow.  If the eigenvalues of $u_1$ are $0$ or $1$, this implies that $dim(Ker(u_1))\geq 2)$ or $dim(Ker(u_1)-Id_{\mathbb{R}^4})\geq 2$,  we can apply  the proposition 5.4 and the proposition 5.5 to deduce that $\pi_1(B(M))$ is nilpotent.

\medskip

{\bf Theorem 5.2.}
{\it Let $(M,P\nabla_M)$ be a $3$-dimensional projective manifold whose developing map is injective, suppose that $dim(Aut(M,P\nabla_M))\geq 2$, then $M$ is homeomorphic to a spherical manifold, $S^2\times S^1$ or a finite cover of $M$ is a torus bundle.}

\medskip

{\bf Proof.}
Suppose that $Aut(M,P\nabla_M)_0$ is not solvable, the proposition 5.2 and the proposition 5.3 imply that either  $M$ is the total space of a bundle over $S^1$ whose fibre is homeomorphic to $T^2$  or $\pi_1(M)$ is solvable. If $Aut(M,P\nabla_M)_0$ is solvable, the theorem 5.1, the propositions 5.4, 5.5 and 5.6 show that $\pi_1(M)$ is solvable. 
In [2], it is shown that a $3$-dimensional closed manifold whose fundamental group is solvable is homeomorphic to a spherical manifold, $S^1\times S^2$, a finite cover of $M$ is a torus bundle,  or $\mathbb{R}P^3\# \mathbb{R}P^3$. Benoist [3] and Goldman and Cooper [8]  have shown that there does not exist a projective structure on $\mathbb{R}P^3\# \mathbb{R}P^3$.

\bigskip

{\bf References.}

\bigskip

1. Albert A. Structure of algebras. Vol. 24. American Math. Soc. 1939.

\smallskip

2. Aschenbrenner, M. Friedl, S. Wilton, H. 3 manifold groups. https://arxiv.org/pdf/1205.0202.pdf

\smallskip

3. Benoist, Y. Nilvariétés projectives. Commentarii Mathematici Helvetici, 1994, vol. 69, no 1, p. 447-473.

\smallskip

4. Benzecri J.P. Sur les vari\'et\'es localement affines et localement projectives.
{\it Bulletin de la S.M.F.} 88 (1960) 229-332.

\smallskip

5. Bredon, G. Bredon, Glen E. Topology and geometry. Vol. 139. Graduate texts in Math.

\smallskip

6. Chatelet, Gilles, Rosenberg, Harold, et Weil, Daniel. A classification of the topological types of $\mathbf {R}^ 2$-actions on closed orientable 3-manifolds. Publications Mathématiques de l'IHÉS, 1974, vol. 43, p. 261-272

\smallskip

7. Conlon, Lawrence. Transversally parallelizable foliations of codimension two. Transactions of the American Mathematical Society, 1974, vol. 194, p. 79-102.

\smallskip

8. Cooper, D.  Goldman, W. A 3–Manifold with no Real Projective Structure. In Annales de la Faculté des sciences de Toulouse: Mathématiques, vol. 24, no. 5, pp. 1219-1238. 2015

\smallskip

9. Goldman, W. Convex real projective structures on compact surfaces.
{\it J. Differential Geometry} 31 (1990) 791-845.

\smallskip

10. Montgomery D. Samelson H. Transformation groups of spheres. {\it Annals of Mathematics} 44 (1943) 454-470.

\smallskip

11. Montgomery, D. Zippin, L. Topological transformation groups. Interscience Tracts in pure and applied mathematics 1955.

\smallskip

12. Sullivan, D;  Thurston,, W. Manifolds with canonical coordinate charts: some examples." Enseign. Math 29 (1983): 15-25.

\smallskip

13. Tsemo, A. Dynamique des vari\'et\'es affines. {\it  Journal of the London Mathematical Society} 63 (20010 469-486.

\smallskip

14. Tsemo, A. D\'ecomposition des vari\'et\'es affines. {\it Bulletin des sciences math\'ematiques} 125 (2001) 71-83.

15. A Tsemo, A. "Linear foliations on affine manifolds." arXiv preprint arXiv:2008.05357 (2020).
\medskip
 
\end{document}